\documentclass[11pt]{article}
\usepackage{latexsym,amsopn,amssymb} 

\newcommand{\qed}{\hfill$\Box$}

\newcommand{\inv}{{\rm inv}}
\newcommand{\D}{{\rm Des}}
\newcommand{\des}{{\rm des}}

\newcommand{\maj}{{\rm maj}}

\renewcommand{\neg}{{\rm neg}}
\DeclareMathOperator{\fm}{flag-major}
\newcommand{\sgn}{{\rm sign}}
\newcommand{\last}{{\rm last}}

\newcommand{\bpi}{{\bar\pi}}

\newcommand{\eps}{{\varepsilon}}
\renewcommand{\ell}{l}

\newtheorem{thm}{Theorem}[section]

\newtheorem{lem}[thm]{Lemma}

\newtheorem{cor}[thm]{Corollary}

\newtheorem{rem}[thm]{Remark}

\def\r#1#2{(#1)_{#2}}
\newcommand{\z}{\zeta}
\def\en#1{\sum_{n=0}^\infty #1 \frac{u^n}{\r qn}}
\def\ent#1{\sum_{n=0}^\infty #1{u^n}/{\r qn}}

\begin{document}

\title{Signed Mahonians}
\bibliographystyle{acm}
\author{Ron M.\ Adin%
\thanks{Department of Mathematics and Statistics, Bar-Ilan University,
Ramat-Gan 52900, Israel. 
Email: {\tt radin@math.biu.ac.il} } 
\and Ira M.\ Gessel
\thanks{Department of Mathematics, Brandeis University,
Waltham,  MA 02454, USA. 
Email: {\tt gessel@brandeis.edu} } 
\and Yuval Roichman%
\thanks{Department of Mathematics and Statistics, Bar-Ilan University,
Ramat-Gan 52900, Israel. 
Email: {\tt yuvalr@math.biu.ac.il} } 
\thanks{The first and the third authors were partially supported by the EC's IHRP Programme, 
within the Research Training Network ``Algebraic Combinatorics in Europe'', 
grant HPRN-CT-2001-00272. The second author was partially supported by NSF Grant DMS-0200596.}}
\date{February 12, 2004}

\maketitle

\begin{abstract}
A classical result of MacMahon gives a simple product formula for 
the generating function of major index over the symmetric group. 
A similar factorial-type product formula for the generating function of 
major index together with sign was given by Gessel and Simion.
Several extensions are given in this paper, including a
recurrence formula, a specialization at roots of unity
and type $B$ analogues.
\end{abstract}

\section{Introduction}

\subsection{Outline}

Enumeration over the symmetric group $S_n$ and related combinatorial objects, taking into account also the {\em sign} of each permutation, was studied by Simion and Schmidt~\cite{SiS} and others
(see, e.g., \cite{Ruskey, White, St.balance, AR321, Reif}).

The polynomial 
$$
\sum_{\pi\in S_n} \sgn(\pi)\, q^{\des(\pi)}
$$
was called the {\em signed Eulerian} by D\'esarmenien and Foata~\cite{DF}.
An elegant formula for signed Eulerians, conjectured by Loday~\cite{Lod}, 
was proved by D\'esarmenien and Foata~\cite{DF} and by Wachs~\cite{Wachs}.
Type $B$ analogues were given by Reiner~\cite{Rein}.

\medskip

MacMahon showed, about a hundred  years ago, that the generating function for major
index over the symmetric group has a simple product formula. The 
{\it signed Mahonian}
will be defined as
$$
\sum_{\pi\in S_n} \sgn(\pi)\, q^{\maj(\pi)}.
$$
An elegant factorial-type product formula for
the signed Mahonians was given by 
Gessel and Simion~\cite[Cor.~2]{Wachs} (Theorem~\ref{t.GS} below).
Various extensions of this theorem are given in
this paper.

\medskip

First, a recurrence for the joint distribution of the inversion number,
major index, and last digit of a permutation is given (Theorem~\ref{t.rec} 
below). It is shown that these parameters give rise to a multiplicative,  
factorial-type formula, if the parameter for inversion number is set equal to
$1$ or to $-1$ (Theorem~\ref{t.gf} below). 

\smallskip

An extension in a different direction gives a factorization of the bivariate 
generating function of major index and inversion number at roots of unity 
(Theorem~\ref{t.g1} below). The proof applies a remarkable identity which follows
from results of Gordon~\cite{Go}, Roselle~\cite{Ros}, and 
Foata-Sch\"utzenberger~\cite{FS}. The identity was independently proved
by Gessel~\cite[Theorem 8.5]{Ge}.

\smallskip

These extensions imply two 
different new proofs of Theorem~\ref{t.GS}.

\medskip

Then Theorem~\ref{t.GS} is extended to the group of signed permutations
$B_n$, where the generating function
of the flag-major index with each of the one-dimensional characters
is shown to have a similar factorial type formula (Theorems~\ref{t.fm1},
\ref{t.fm2} and~\ref{t.fm3} below). 

\medskip

These results yield explicit simple generating functions for the (flag) 
major index on subgroups of index 2 of $S_n$ and $B_n$, such as 
the alternating groups and the Weyl groups of type $D$.
See Section~\ref{s.subgps}.

\bigskip

The rest of the paper is organized as follows.
Necessary background and statements of main results are given
in the rest of this section.
In Section~\ref{s.rec}, a multivariate recurrence formula 
for length, major index and last digit is proved (Theorem~\ref{t.rec}). 
Then, in Section~\ref{s.gf}, this formula is applied to prove 
a new extension (Theorem~\ref{t.gf}) of the Gessel-Simion Theorem.
A second proof of the Gessel-Simion Theorem, via specialization at 
roots of unity, is given in Section~\ref{s.root}. 
The type $B$ analogue (Theorem~\ref{t.fm1}) is proved in Section~\ref{s.bsup}.
The distribution of the (flag) major index on 
index~2 subgroups is then deduced in Sections~\ref{s.char} and~\ref{s.subgps}.

\subsection{Background}

The Coxeter generators  $\{\, s_i=(i,i+1)\mid\ 1\le i\le n-1\,\}$
of $S_n$ give rise to various combinatorial statistics.
For $\pi\in S_n$
let  the {\it length}, $\ell(\pi)$, be the standard length of $\pi$ with respect
to these generators, which is the same as the number of inversions of $\pi$.
This notion is defined similarly for other Coxeter groups.
The generating function of length in a Coxeter group $W$
is called the {\it Poincar\'e polynomial} of $W$~\cite[Ch. 3]{Hu}.


For a positive integer $n$ define 
$$
[n]_q\ :=\ \ {1-q^n\over 1-q}.
$$
Then 

\begin{thm}\label{t.po}{\rm\cite[\S 3.15]{Hu}}
$$
\sum\limits_{\pi\in S_n} q^{\ell(\pi)}=
[1]_q [2]_q \cdots [n]_q. 
$$
\end{thm}

Another statistic on $S_n$, which has a Coxeter group interpretation, 
is the descent number.
Given a permutation $\pi$ in the symmetric group 
$S_n$, the {\it descent set} of $\pi$ is 
$$
\D (\pi):=\{\,i \mid \ell(\pi)>\ell(\pi s_i)\,\}
=\{\,i\mid \pi (i) > \pi (i+1)\,\}
$$
\noindent
and the corresponding {\it descent number} is $\des(\pi):=| \D(\pi) |$.
The {\it major index} of $\pi$ is 
the following weighted enumeration of the descents
$$ 
\maj (\pi) := \sum\limits_{i\in \D (\pi)} i.
$$
A well-known classical result 
asserts that
the length function and major index of a permutation are equidistributed
over the symmetric group $S_{n}$. 

\begin{thm}\label{t.mm}{\rm(MacMahon~\cite{MM})}
$$
\sum\limits_{\pi\in S_n} q^{\ell(\pi)}=
\sum\limits_{\pi\in S_n} q^{\maj(\pi)}=
[1]_q [2]_q \cdots [n]_q.
$$
\end{thm}


\bigskip

A similar simple factorial-type product formula for
the signed Mahonians was given by 
Gessel and Simion~\cite[Cor.~2]{Wachs}. 

The {\it sign} of an element $w$ in a Coxeter group $W$ is 
$$
\sgn(w)\ :=\ \ (-1)^{\ell(w)}.
$$

\begin{thm}\label{t.GS}{\rm(The Gessel-Simion Theorem)}
$$
\sum\limits_{\pi\in S_n} \sgn(\pi) q^{\maj(\pi)} 
= [1]_q [2]_{-q} [3]_q [4]_{-q}\cdots [n]_{(-1)^{n-1}q}.
$$
\end{thm}

\bigskip

Recall that $B_{n}$ denotes the group of all bijections $\sigma$ of the set 
$[-n,n] \setminus \{ 0 \}$ onto itself such that 
$$
\sigma (-a)=-\sigma (a)
$$
for all $a \in [-n,n] \setminus \{0\}$, with composition as the group operation.
This group is usually known as the group of ``signed permutations'' on $[n]$, 
or as the {\em hyperoctahedral group} of rank $n$,
or as the classical Weyl group of type $B$ and rank $n$. 

\bigskip

It is well known (see, e.g., \cite[Proposition 8.1.3]{BB}) that $B_{n}$
is a Coxeter group with respect to the generating set $\{ s_{0},s_{1},s_{2},
\ldots , s_{n-1} \}$, where
\[ s_{0} := [-1,2, \ldots n ] \]
and
\[ s_{i} := [1,2,\ldots ,i-1,i+1, i,i+2, \ldots n ] \]
for $i=1, \ldots , n-1$. 
Let $\ell(\sigma)$ be the standard length of $\sigma\in B_n$
with respect to its Coxeter generators.

\begin{thm}\label{t.poB}{\rm\cite[\S 3.15]{Hu}}
$$
\sum\limits_{\pi\in B_n} q^{\ell(\pi)}=
[2]_q [4]_q \cdots [2n]_q. 
$$
\end{thm}

\medskip

Despite the fact that an increasing number of enumerative results 
of this nature have been generalized to the hyperoctahedral group $B_n$ 
(see, e.g., \cite{Bre94, FH, Rei93c, Rei95, 
Sta76}) 
and that several ``major index'' statistics have been introduced 
and studied for $B_n$~\cite{CF94, CF95a, CF95b, FK, Rei93a, 
Rei93b, Ste94} no generalization of MacMahon's result to $B_n$ has been found until a new statistic, the
{\it flag major index}, was introduced.

\medskip

The {\it flag-major} of $\sigma\in B_n$ is defined as
$$
\fm(\sigma)\ :=\ \  2\cdot \maj(\sigma)+\neg(\sigma)
$$
where 
$$
\neg(\sigma)\ :=\ \ \#\{\,1\le i\le n \mid \sigma(i)<0\,\}
$$
and $\maj(\sigma)$ is the major index of the sequence
$(\sigma(1),\dots,\sigma(n))$, with respect to the order
$$
-1<\cdots<-n<1<\cdots< n .
$$
A type $B$ analogue of Theorem~\ref{t.mm}
was given in~\cite{AR}.

\begin{thm}\label{t.ar}{\rm\cite{AR}}
$$
\sum\limits_{\pi\in B_n} q^{\ell(\pi)}=
\sum\limits_{\pi\in B_n} q^{\fm(\pi)}=
[2]_q [4]_q \cdots [2n]_q.
$$
\end{thm}

For a unified definition of the classical major index and the flag-major index 
as a length of a distinguished canonical expression see~\cite{AR}. 
The flag-major index has many combinatorial and algebraic properties
which are shared with the classical major index on $S_n$~\cite{AR, ABR, HLR, ABR2, ABR4, BC1, CG}.
In this paper we will give a type $B$ analogue of the Gessel-Simion Theorem (Theorem~\ref{t.fm1} below),
as well as other new extensions of this theorem.

\subsection{Main Results}

We find a recurrence (theorem~\ref{t.rec} below)
for the joint distribution of length, major index, and 
last digit, which leads to the following result.
Let
$$
\last(\pi)\ :=\ \ \pi(n)-1.
$$
Then

\begin{thm}\label{t.main3}{\rm(see Theorem~\ref{t.gf} below)}\\
For $\eps = \pm1$,
$$
\sum_{\pi\in S_n} \eps^{\ell(\pi)} q^{\maj(\pi)} z^{\last(\pi)}
= [1]_q \cdot [2]_{\eps q} \cdot [3]_q \cdot [4]_{\eps q} \cdots 
  [n-1]_{\pm q} \cdot [n]_{\pm\eps q/z} \cdot z^{n-1}.
$$
\end{thm}

This theorem shows that the distribution of (signed) major index over 
permutations with prescribed last digit is essentially independent of 
this digit (Corollary~\ref{t.last-indep}). 
Letting $\eps=-1$ and $z=1$ gives Theorem~\ref{t.GS}.

\bigskip

A second new proof of Theorem~\ref{t.GS} uses a known identity
(Theorem~\ref{t.ros} below) involving the 
generating function for length and major index. 
This also leads to 
a factorization at roots of unity other than $\pm1$.

Let
$$
A_n(t,q)\ :=\ \ \sum\limits_{\pi\in S_n} t^{\ell(\pi)} q^{\maj(\pi)}.
$$
For a positive integer $n$ define
$$
(q)_n\ :=\ \ (1-q)(1-q^2)\cdots (1-q^n).
$$

\begin{thm}\label{t.ges}{\rm(see Theorem~\ref{t.g1} below)}\\
Let $n$ and $m$ be positive integers.
Let $\z$ be a primitive $m$th root of unity  and assume that
$n=mk+i$ with $0\le i<m$. Then
$$
A_n(\z,q)= A_i(\z,q){(q)_n\over (q)_i (1-q^m)^k}.
$$
\end{thm}

The case $m=2$ gives Theorem~\ref{t.GS}.

\bigskip

A type $B$ analogue of Theorem~\ref{t.GS} is :

\begin{thm}\label{t.main2}{\rm(see Theorem~\ref{t.fm1} below)}
$$
\sum\limits_{\pi\in B_n} \sgn(\pi)\cdot q^{\fm(\pi)} 
=  [2]_{- q} [4]_{ q} \cdots [2n]_{(-1)^{n} q}.
$$
\end{thm} 

\medskip

Explicit generating functions of the major index and
flag major index
on distinguished subgroups follow from Theorems~\ref{t.GS}
and~\ref{t.main2}. See Corollaries~\ref{t.subgps1} and~\ref{t.subgps2} below.

\section{A Recurrence Formula}\label{s.rec}

Let $S_n$ be the symmetric group.
For $\pi \in S_n$ define the following statistics:
\begin{eqnarray*}
\inv(\pi) &:=& \mbox{\rm inversion number of $\pi$}\\ 
 &(=& \mbox{\rm length of $\pi$ w.r.t. the usual Coxeter generators of $S_n$})\\
\maj(\pi) &:=& \mbox{\rm major index of $\pi$} = \sum \{\,1\le i\le n-1\mid \pi(i) > \pi(i+1)\,\}\\
\last(\pi) &:=& \pi(n)-1, \mbox{\rm\ one less than the last digit in $\pi$}
\end{eqnarray*}
Define the multivariate generating function
\begin{equation}\label{e.gf}
f_n(x,y,z) := \sum_{\pi\in S_n} x^{\inv(\pi)} y^{\maj(\pi)} z^{\last(\pi)}.
\end{equation}
\begin{thm}\label{t.rec}{\bf (Recurrence Formula)}
$$
f_1(x,y,z) = 1
$$
and, for $n\ge 2$,
\begin{eqnarray*}
(x-z)f_n(x,y,z) &=& (x^{n}y^{n-1} - z^{n}) \cdot f_{n-1}(x,y,1)\\
                       && \quad+\; x^{n-1} (1 - y^{n-1}) z
\cdot f_{n-1}(x,y,z/x).
\end{eqnarray*}
\end{thm}

\noindent
{\bf Proof.}
The case $n=1$ is clear. Assume $n\ge 2$.
\par\noindent 
Given a permutation
$$
\pi = (\pi(1),\ldots,\pi(n-1))\in S_{n-1},
$$
append $k$ ($1\le k\le n$) as the $n$th digit, while adding $1$ to each 
existing digit between $k$ and $n-1$, to get a permutation
$$
\bpi = (\bpi(1),\ldots,\bpi(n-1),k)\in S_n
$$
where, for $1\le i\le n-1$,
$$
\bpi(i) = \cases{%
\pi(i), &if $\pi(i) < k$;\cr
\pi(i)+1, &otherwise.}
$$
The new statistics for $\bpi$ are:
\begin{eqnarray*}
\inv(\bpi) &=& \inv(\pi) + (n-k)\\
\maj(\bpi) &=& \cases{\maj(\pi),& if $k > \pi(n-1)$;\cr \maj(\pi)+(n-1),& otherwise.}\\
\last(\bpi) &=& k-1
\end{eqnarray*}

We can therefore compute
\begin{eqnarray*}
f_n &=& f_n(x,y,z)\\ 
&=& \sum_{\pi\in S_{n-1}} \sum_{k=1}^{n} x^{\inv(\bpi)} y^{\maj(\bpi)} z^{\last(\bpi)}\\
&=& \sum_{\pi\in S_{n-1}} x^{\inv(\pi)+n-1} y^{\maj(\pi)} \\ 
& &\quad \times \left[y^{n-1} \sum_{k=1}^{\last(\pi)+1}
x^{1-k} z^{k-1} 
                               + \sum_{k=\last(\pi)+2}^{n} x^{1-k} z^{k-1}\right]\\
&=& (1-z/x)^{-1}\sum_{\pi\in S_{n-1}} x^{\inv(\pi)+n-1} y^{\maj(\pi)} \\
& &\quad\times \left[y^{n-1}\left(1 -
(z/x)^{\last(\pi)+1}\right) +
\left((z/x)^{\last(\pi)+1} - (z/x)^{n}\right)\right]\\ &=&
(1-z/x)^{-1} [ (x^{n-1} y^{n-1} - x^{-1} z^{n})
f_{n-1}(x,y,1)\\ & &\quad+ \left. x^{n-2} (1 - y^{n-1}) z
f_{n-1}(x,y,z/x) \right..\\
\end{eqnarray*}
Multiplying both sides by $x-z$ gives the claimed recurrence.

\qed

\section{A Multiplicative Generating Function}\label{s.gf}

In general, the generating function from the previous section is a complicated polynomial of its variables. However, assuming in addition that $x^2 = 1$ leads to
surprisingly simple  results.
\begin{cor}\label{t.fsmall}
The first few values of $f_n$, assuming $x = \eps = \pm1$, are:
\begin{eqnarray*}
f_1(\eps,q,z) &=& 1\\
f_2(\eps,q,z) &=& z + \eps q\\
f_3(\eps,q,z) &=& (1 + \eps q)(z^2 + qz + q^2)\\
f_4(\eps,q,z) &=& (1 + \eps q)(1 + q + q^2)(z^3 + \eps qz^2 + q^2z + \eps q^3)
\end{eqnarray*}
\end{cor}

The case $\eps = z = 1$ is a well-known result of MacMahon.

\begin{thm}\label{t.gf}
For $\eps = \pm1$,
\begin{eqnarray*}
\sum_{\pi\in S_n} \eps^{\inv(\pi)} q^{\maj(\pi)} z^{\last(\pi)}
&=& \left( \prod_{i=1}^{n-1} [i]_{\eps^{i-1} q} \right) \cdot [n]_{\eps^{n-1} q/z} \cdot z^{n-1}\\
&=& [1]_q [2]_{\eps q} [3]_q [4]_{\eps q} \cdots [n-1]_{\pm q} \cdot [n]_{\pm\eps q/z} z^{n-1}.
\end{eqnarray*}
\end{thm}

\noindent
{\bf Proof.}
By induction on $n$. 
By Corollary~\ref{t.fsmall}, the claim is true for $n=1$ (as well as for $n=2,3,4$). Assume now that the claim holds for $n-1$, where $n\ge 2$.
Thus
$$
f_{n-1}(\eps,q,z) =
\left( \prod_{i=1}^{n-2} [i]_{\eps^{i-1} q} \right) \cdot [n-1]_{\eps^{n-2} q/z} \cdot z^{n-2}.
$$
Substituting in the recurrence formula of Theorem~\ref{t.rec} and eliminating the factor
$$
\left( \prod_{i=1}^{n-2} [i]_{\eps^{i-1} q} \right),
$$
it remains to show that
\begin{eqnarray*}
& & (\eps - z) [n-1]_{\eps^{n-2} q} [n]_{\eps^{n-1}q/z} \cdot z^{n-1}\\
&=& (\eps^{n}q^{n-1} - z^{n}) [n-1]_{\eps^{n-2} q}
+ \eps^{n-1} (1 - q^{n-1}) z [n-1]_{\eps^{n-1} q/z} \cdot (z/\eps)^{n-2}.
\end{eqnarray*}
Using the definition of $[k]_q$,
this is equivalent to
\begin{eqnarray*}
& & \frac{(\eps - z) (1-(\eps^{n-2} q)^{n-1}) (z^{n}-(\eps^{n-1} q)^{n})}{%
    (1-\eps^{n-2} q) (z-\eps^{n-1} q)}\\
&=& \frac{(\eps^{n}q^{n-1} - z^{n}) (1-(\eps^{n-2} q)^{n-1})}{1-\eps^{n-2} q}
    + \frac{\eps z (1 - q^{n-1}) (z^{n-1}-(\eps^{n-1} q)^{n-1})}{z-\eps^{n-1} q}.
\end{eqnarray*}
Clearing denominators and using the fact that $(n-2)(n-1)$ is even, we can transform this equation into
\begin{eqnarray*}
(\eps - z) (1-q^{n-1}) (z^{n}-q^{n})
&=& (\eps^{n}q^{n-1} - z^{n}) (1-q^{n-1}) (z-\eps^{n-1} q)\\
&+& \eps z (1 - q^{n-1}) (z^{n-1}-\eps^{n-1}q^{n-1}) (1-\eps^{n-2} q).
\end{eqnarray*}
Dividing by $(1-q^{n-1})$ one gets
\begin{eqnarray*}
(\eps - z) (z^{n}-q^{n})=
(\eps^{n}q^{n-1} - z^{n}) (z-\eps^{n-1} q)+
\eps z (z^{n-1}-\eps^{n-1}q^{n-1}) (1-\eps^{n} q),
\end{eqnarray*}
completing the proof.

\qed

Letting $z=1$, one gets
\begin{cor}\label{t.maj3}
\begin{eqnarray*}
\sum_{\pi\in S_n} q^{\maj(\pi)} 
&=& [n]_{q}! := [1]_{q} [2]_{q} \cdots [n]_{q}\\
\sum_{\pi\in S_n} \sgn(\pi)\, q^{\maj(\pi)} 
&=& [n]_{\pm q}! := [1]_{q} [2]_{-q} [3]_{q} [4]_{-q} \cdots [n]_{(-1)^{n-1} q}
\end{eqnarray*}
\end{cor}

The first formula is a classical result of MacMahon~\cite{MM}, 
and the second was first proved by Gessel and Simion~\cite[Cor. 2]{Wachs}.

\begin{cor}\label{t.last-indep}
The distributions of $\maj$ and of $\maj$ with $\sgn$ over all 
permutations with a prescribed last digit are essentially independent of 
this digit, namely: if
$$
S_n(k) := \{\,\pi\in S_n \mid \pi(n) = k\,\}\qquad(1\le k\le n)
$$
then, for $\eps = \pm1$,
\begin{eqnarray*}
\sum_{\pi\in S_n(k)} \eps^{\inv(\pi)} q^{\maj(\pi)} 
&=& f_{n-1}(\eps,q,1) \cdot (\eps^{n-1} q)^{n-k}\\
&=& \left( \prod_{i=1}^{n-1} [i]_{\eps^{i-1} q} \right) \cdot (\eps^{n-1} q)^{n-k}.
\end{eqnarray*}
\end{cor}

\newpage

\section{Specialization at Roots of Unity}\label{s.root}

A proof of Theorem~\ref{t.ges} is given in this section. 

\bigskip

Suppose that we have a sequence $f_0(q), f_1(q),\dots$ of
polynomials in
$q$ defined by  the Eulerian generating function 
\begin{equation} \label{e:egf} 
F(u;q)=\en {f_n(q)},
\end{equation}
where $\r qn := (1-q)(1-q^2)\cdots (1-q^n)$. 
We would like to study the values of $f_n(q)$ at a root of unity. We cannot 
simply evaluate (\ref{e:egf}) at a root of unity, since this would make
denominators vanish. Instead we take a less direct approach.

Fix a positive integer $m$, and 
let $\phi_m(q)$ be the cyclotomic polynomial of order $m$ in $q$
(whose roots are all the primitive $m$th roots of unity).
If $f(q)$ and $g(q)$ are polynomials in $q$ with rational coefficients 
and $\z$ is a primitive $m$th root of unity, then 
$f(q)\equiv g(q)\pmod {\phi_m(q)}$ if and only if $f(\z)=g(\z)$, 
since $\phi_m(q)$ is irreducible over the rationals and $\phi_m(\z)=0$.

Given two Eulerian generating functions 
$F(u;q)=\ent{f_n(q)}$ and $G(u;q)=\ent{g_n(q)}$, 
by $F(u;q)\equiv G(u;q)$ we mean that
$f_n(q)\equiv g_n(q)\pmod{\phi_m(q)}$ for all $n$. 
Henceforth we take all congruences to be modulo $\phi_m(q)$.

The basic facts about these congruences are contained in the following lemma:
\begin{lem}\label{t.lemma1}
 Let $u_i := u^i/\r qi$.
\begin{enumerate}
\item[(i)] If $0\le i,j < m$ and $i+j\ge m$ then $u_iu_j\equiv 0$.
\item[(ii)] If $0\le i < m$ then 
$$u_{mk+i}\equiv {u_m^k\over k!}{u_i}.$$
\end{enumerate}
\end{lem}

\noindent{\bf Proof.} 
Let $\z$ be a primitive $m$th root of unity. For (i), we have
$$
u_iu_j = {u^{i+j}\over (q)_i (q)_j} = {(q)_{i+j}\over (q)_i(q)_j} u_{i+j}.
$$
The quotient in the right-hand-side is a polynomial in $q$ 
(a $q$-binomial coefficient, see below).
Since $(q)_{i+j}$ vanishes for $q=\z$ but $(q)_i(q)_j$ does not, (i) follows.

For (ii), we have 
$${u_m^k\over k!}u_i={u^{mk+i}\over (q)_m^k k!\,(q)_i}
  = {(q)_{mk+i}\over (q)_m^k k!\,(q)_i} u_{mk+i},
$$
so it suffices to show that 
$$
\left.{(q)_{mk+i}\over (q)_m^k k!\,(q)_i} \right|_{q=\z}=1.
$$
To prove this, we show that 
$$
\left.{(q)_{mk+i}\over (q)_{mk}(q)_i} \right|_{q=\z}=1
$$
and that 
$$
\left. {(q)_{mk}\over (q)_m^k}\right|_{q=\z}=k!.
$$

For the first equality, we have
$$
{(q)_{mk+i}\over (q)_{mk}(q)_i}={1-q^{mk+1}\over 1-q}
{1-q^{mk+2}\over 1-q^2}
\cdots {1-q^{mk+i}\over 1-q^i}.
$$
Since $\z^{mk+j}=\z^{j}\ne1$ for $j=1,2,\dots, i$, the equality follows.

For the second equality, let us write
$$
(q)_{mk}=\prod_{1\le l\le mk\atop m\nmid l}
(1-q^l)\cdot\prod_{j=1}^k (1-q^{mj}),
$$
so
$$
{(q)_{mk}\over (q)_m^k}={\displaystyle\prod_{1\le l\le mk,\, m\nmid l} 
(1-q^l)\over (q)_{m-1}^k} \cdot\prod_{j=1}^k{1-q^{mj}\over 1-q^m}.
$$

We may evaluate the first factor on the right at $q=\z$ by simply setting 
$q=\z$, since neither the numerator nor the denominator vanishes, and we see 
easily that this factor becomes 1. Writing the second factor as 
$$
\prod_{j=1}^{k} (1+q^m+\cdots +q^{m(j-1)}),
$$ 
we see that setting $q=\z$ in it yields $k!$.
\qed

Now recall that the {$q$-binomial coefficient} ${n\brack k}_q$ is 
the polynomial in $q$ defined by 
$$
{n\brack k}_q = {(q)_n\over (q)_k (q)_{n-k}}
$$
for $0\le k\le n$, with ${n\brack k}_q=0$ for $n<k$.
As a consequence of Lemma~\ref{t.lemma1} we obtain a frequently rediscovered 
result of Gloria Olive \cite[(1.2.4)]{Olive} about the evaluation of 
$q$-binomial coefficients at roots of unity:

\begin{cor}\label{t.olive}
Let $m$ be a positive integer and let $\z$ be a primitive $m$th root of unity. 
Let $a_1$, $a_2$, $b_1$, and $b_2$ be nonnegative integers with $0\le b_1,b_2< m$. 
Then
$$
{(ma_1+b_1)+(ma_2+b_2)\brack ma_1+b_1}_\z 
= {a_1+a_2\choose a_1}{b_1+b_2\brack b_1}_\z.
$$
\end{cor}

\noindent{\bf Proof.} With the notation of Lemma~\ref{t.lemma1} we have

\begin{eqnarray*}
{(ma_1+b_1)+(ma_2+b_2)\brack ma_1+b_1}_q u_{(ma_1+b_1)+(ma_2+b_2)}
& = &{u^{ma_1+b_1}\over(q)_{ma_1+b_1}} {u^{ma_2+b_2}\over(q)_{ma_2+b_2}} \\[5pt]
& = & u_{ma_1+b_1} u_{ma_2+b_2}.
\end{eqnarray*}
By Lemma~\ref{t.lemma1}(ii) this is congruent modulo $\phi_m(q)$ to 
$$
{u_m^{a_1}\over a_1!}u_{b_1} {u_m^{a_2}\over a_2!}u_{b_2}.
$$
If $b_1+b_2\ge m$ then, by Lemma~\ref{t.lemma1}(i), this is congruent to 0. 
Otherwise we have, by Lemma~\ref{t.lemma1}(ii),
\begin{eqnarray*}
{u_m^{a_1}\over a_1!}u_{b_1} {u_m^{a_2}\over a_2!}u_{b_2}
&=& {a_1+a_2\choose a_1} {u_m^{a_1+a_2}\over (a_1+a_2)!} \cdot 
    {b_1+b_2\brack b_1}_q u_{b_1+b_2}\\[3pt]
&\equiv& {a_1+a_2\choose a_1}{b_1+b_2\brack b_1}_q u_{m(a_1+a_2)+(b_1+b_2)},
\end{eqnarray*}
and the result follows.

\qed

Our proof of Theorem \ref{t.ges} is based on the  generating function for the
bivariate distribution of length and major index:

\begin{thm}\label{t.ros}
Let the polynomials $A_n(q,r)$ be defined by
\begin{equation}\label{e:ip} 
A(u;q):=\prod_{i,j=0}^\infty {1\over 1-q^ir^ju}=
\sum_{n=0}^\infty {A_n(q,r)\over (q)_n(r)_n}u^n.
\end{equation} 
Then
$$
A_n(q,r)\ =\ \ \sum\limits_{\pi\in S_n} q^{\ell(\pi)} r^{\maj(\pi)}.
$$
\end{thm}

\noindent{\bf Historical Note:}
Theorem \ref{t.ros} was first proved by Gessel~\cite[Theorem 8.5]{Ge}. 
(For a refinement that also includes the number of descents, see~\cite{GG}.)
Basil Gordon~\cite{Go} had earlier given a combinatorial interpretation
to the coefficients of $A_n(q,r)$, but did not describe it very explicitly. 
(In fact, he considered the generalization 
$\prod_{i,j,\dots, k=0}^\infty ( 1-q^ir^j\cdots s^k u)^{-1}$.)
D.~P.~Roselle~\cite{Ros} 
explained Gordon's combinatorial interpretation more explicitly.
His result is equivalent to
$$
A_n(q,r)=\sum_{\pi\in S_n} q^{\maj(\pi^{-1})}r^{\maj (\pi)}.
$$ 
Then D.~Foata and M.-P.~Sch\"utzenberger~\cite{FS}  
gave a bijective proof that 
$$\sum_{\pi\in S_n} q^{\maj (\pi^{-1})}r^{\maj(\pi)}
  =\sum_{\pi\in S_n} q^{l(\pi)}r^{\maj(\pi)},
$$ 
which, together with the result of Gordon and Roselle, implies 
Theorem~\ref{t.ros}.

\begin{thm}\label{t.g1} 
Let $n$ and $m$ be positive integers.
Let $\z$ be a primitive $m$th root of unity, and assume that
$n=mk+i$ with $0\le i<m$. Then
$$
A_n(\z,r)= A_i(\z,r){(r)_n\over (r)_i (1-r^m)^k}.
$$
\end{thm}

\noindent{\bf Proof.}
To find a congruence modulo $\phi_m(q)$ for the polynomials
$A_n(q,r)$, think of (\ref{e:ip}) as an Eulerian generating function 
in which the coefficient of $u^n/\r qn$ is $A_n(q,r)/\r rn$. By taking 
logarithms and exponentiating, we see that 
\begin{eqnarray*}
A(u;q) 
&=& \prod_{i,j=0}^{\infty} {1 \over 1 - q^i r^j u}
\;=\; \exp\biggl(-\sum_{i,j=0}^{\infty} \ln\,(1 - q^i r^j u)\biggr)\\
&=& \exp\biggl(\,\sum_{i,j=0}^{\infty} \sum_{t=1}^{\infty} 
    {(q^i r^j u)^t \over t}\biggr)
\;=\; \exp\biggl(\sum_{t=1}^{\infty} {u^t \over t(1-q^t)(1-r^t)}\biggr).
\end{eqnarray*}
Now 
$$
\sum_{t=1}^\infty {u^t \over t(1-q^t)(1-r^t)} 
= \sum_{t=1}^\infty {\r q{t-1}\over t(1-r^t)}{u^t\over \r qt}
\equiv \sum_{t=1}^m {\r q{t-1}\over t(1-r^t)}{u^t\over \r qt},
$$ 
so
$$
A(u;q)
\equiv \exp\biggl(\sum_{t=1}^{m-1} {\r q{t-1} \over t(1-r^t)} u_t\biggr)
\cdot \exp\biggl({\r q{m-1} \over m(1-r^m)} u_m\biggr).
$$
Using Lemma~\ref{t.lemma1}(i) we see that 
$$
\exp\biggl(\sum_{t=1}^{m-1} {\r q{t-1} \over t(1-r^t)} u_t\biggr)
\equiv \sum_{i=0}^{m-1} B_i(q,r) u_i,
$$ 
where $B_i(q,r)$ are polynomials in $q$ whose coefficients are 
rational functions of $r$.

Now let $\z$ be a primitive $m$th root of unity. 
Setting $x= 1$ in 
$$
(1-\z x)\cdots(1-\z^{m-1}x) = (1-x^{m})/(1-x) = 1+x+\cdots +x^{m-1}$$
we see that 
$$
(1-\z)\cdots (1-\z^{m-1}) = m.
$$ 
Thus $\r q{m-1}\equiv m$, so with the terminology of Lemma~\ref{t.lemma1} 
we have
$$
{\r q{m-1} \over m(1-r^m)} u_m \equiv {u_m \over 1-r^m} 
$$ 
and 
$$
\exp\biggl({\r q{m-1}\over m(1-r^m)} u_m\biggr)
\equiv \exp\biggl({u_m\over 1-r^m}\biggr) 
= \sum_{k=0}^\infty {u_m^k\over
k!\,(1-r^m)^k}.
$$ 
It follows that
\begin{eqnarray*}
\sum_{n=0}^\infty {A_n(q,r) \over (r)_n}{u^n \over (q)_n}
&\equiv& \sum_{i=0}^{m-1} \sum_{k=0}^\infty 
         B_i(q,r) {u_iu_m^k\over k!\,(1-r^m)^k }\\
&\equiv& \sum_{i=0}^{m-1} \sum_{k=0}^\infty 
         {B_i(q,r) \over (1-r^m)^k} {u^{mk+i} \over \r q{mk+i}},
\end{eqnarray*}  
by Lemma~\ref{t.lemma1}(ii). Thus, if $n=mk+i$ with $0\le i<m$, then 
$$
{A_n(q,r)\over (r)_n} \equiv {B_i(q,r) \over (1-r^m)^k}
$$
or equivalently
$$
{A_n(\z,r)\over (r)_n} = {B_i(\z,r) \over (1-r^m)^k}.
$$
Letting $k=0$ (so that $n=i$) we get
$$
B_i(\z,r) = {A_i(\z,r) \over (r)_i}\qquad(0\le i<m)
$$
and the result follows.

\qed

\medskip

\noindent{\bf Second Proof of Theorem~\ref{t.GS}.}
Take $m=2$ in Theorem~\ref{t.g1} and simplify.

\qed


For some other results involving the evaluation of $A_n(q,r)$ at roots of unity, see 
\cite{BMS} and \cite{Go}. 
  
\section{A Signed Mahonian for $B_n$}\label{s.bsup}

Let $B_n$ be the hyperoctahedral group.
The {\it flag-major} of $\sigma\in B_n$ is defined as
$$
\fm(\sigma):= 2\,\maj(\sigma)+\neg(\sigma),
$$
where 
$$
\neg(\sigma):=\#\{\,i\mid \sigma(i)<0\,\}
$$
and $\maj(\sigma)$ is the major index of the sequence
$(\sigma(1),\dots,\sigma(n))$, with respect to the order
$$
-1<\cdots<-n<1<\cdots< n .
$$
Recall that for every $\sigma\in B_n$ we define 
$$
\sgn(\sigma)=(-1)^{\ell(\sigma)},
$$
where the length $\ell$ 
(here and throughout this section) is taken with respect to the Coxeter generators of $B_n$.

\begin{thm}\label{t.fm1}
$$
\sum\limits_{\sigma\in B_n} \sgn(\sigma) q^{\fm(\sigma)}=
[2]_{-q} [4]_{q} \cdots [2n]_{(-1)^n q}.
$$
\end{thm}

\begin{rem}\label{r.fm2}
The above order appeared in~\cite{AR}.
In \cite{ABR} we considered another natural order :
$$
-n<\cdots<-1<1<\cdots<n .
$$
The distribution of flag-major is the the same for both orders, 
but the joint distribution of flag-major and length is different, 
and Theorem \ref{t.fm1} does not hold for {\rm flag-major} defined 
with respect to the latter order.
It was shown in \cite{AR} that {\rm flag-major} defined with respect to 
the first order satisfies some further remarkable properties 
(e.g., it is the length of a certain decomposition of the permutation).
These properties do not hold for the second order.
\end{rem}

\noindent{\bf Proof.}
We use the decomposition 
$$
B_n=U_n\cdot S_n,
$$
where 
$$
U_n:=\{\,\tau\in B_n\mid \tau(1)<\cdots<\tau(n)\,\}
$$
with respect to the order
$$
-1<\cdots<-n<1<\cdots< n ,
$$
and
$$
S_n:=\{\,\pi\in B_n\mid \neg(\pi)=0\,\}.
$$
This decomposition appeared in \cite{ABR} (where it was taken with respect
to the other order).

\medskip

Note that every $\sigma\in B_n$ has a unique decomposition
$\sigma=\tau\pi$, $\tau\in U_n, \pi\in S_n$. Then, by definition,
$$
\fm(\sigma)=2\cdot\maj(\pi)+ \neg(\tau).
$$ 
Thus
\begin{eqnarray*}
\sum\limits_{\sigma\in B_n} \sgn(\sigma) q^{\fm(\sigma)}&=&
\sum\limits_{\tau\in U_n,\,\pi\in S_n} \sgn(\tau\pi) q^{2\cdot
\maj(\pi)+\neg(\tau)}\cr
&=& \sum\limits_{\tau\in U_n} \sgn(\tau) q^{\neg(\tau)}
\cdot \sum\limits_{\pi\in S_n} \sgn(\pi) q^{2\cdot \maj(\pi)}.
\end{eqnarray*}
By Corollary~\ref{t.maj3}, the second factor is equal to
$$
\sum\limits_{\pi\in S_n} \sgn(\pi) q^{2\cdot \maj(\pi)}=
[1]_{q^2}[2]_{-q^2}\cdots [n]_{\pm q^2}.
$$
We shall compute the first factor. 
Define
$$
U_n(k):=\{\,\tau\in U_n\mid \neg(\tau)=k\,\} \qquad (0\le k\le n).
$$
Then
$$
\sum\limits_{\tau\in U_n} \sgn(\tau) q^{\neg(\tau)}=
\sum\limits_{k=0}^n  \sum\limits_{\tau\in U_n(k)} \sgn(\tau)\cdot q^k =
 \sum\limits_{k=0}^n q^k\sum\limits_{\tau\in U_n(k)} 
(-1)^{\ell(\tau)}.
$$

Recall from \cite[Proposition 3.1 and Corollary 3.2]{Bre94}
 \cite[Propositions 8.1.1 and 
8.1.2]{BB} 
that
for every $\sigma\in B_n$
$$
\ell(\sigma)=\inv(\sigma)+\sum\limits_{\{\,1\le i\le n\mid\ \sigma(i)<0\,\}}|\sigma(i)|,
$$
where $\inv(\sigma)$ is taken with respect to the order
$$
-n<\cdots<-1<1<\cdots< n .
$$
Now $U_n$ consists of all elements whose entries are increasing with respect to the
order
$-1<\cdots<-n<1<\cdots< n$. Thus for every $\tau\in U_n(k)$
$$
\inv(\tau)={k\choose 2}
$$
and
$$
\ell(\tau)={k\choose 2}+\sum_{i=1}^k |\tau(i)|.
$$
It follows that
\begin{eqnarray*}
\sum\limits_{\tau\in U_n(k)} (-1)^{\ell(\tau)} & = &
\sum\limits_{\tau\in U_n(k)} (-1)^{{k\choose 2}+\sum_{i=1}^k |\tau(i)|}\\
&=&
(-1)^{{k\choose 2}}
\sum\limits_{1\le i_1<\cdots<i_k\le n} (-1)^{\sum_{j=1}^k i_j}.
\end{eqnarray*}
From the $q$-binomial theorem 
$$\prod_{i=1}^{n}(1+q^ix)=\sum_{k=0}^n q^{k+1\choose 2}{n\brack k}_q x^k,$$
it follows that 
$$
\sum\limits_{1\le i_1<\cdots<i_k\le n} q^{\sum_{j=1}^k i_j}=q^{k+1\choose 2}{n \brack
k}_q.
$$
We deduce that
$$
\sum\limits_{\tau\in U_n(k)} \sgn(\tau)=
(-1)^{{k\choose 2}}(-1)^{k+1\choose 2}{n \brack k}_{-1}=
(-1)^k {n \brack k}_{-1},
$$
so
$$
\sum\limits_{\tau\in U_n} \sgn(\tau) q^{\neg(\tau)}=
\sum\limits_{k=0}^n q^k\sum\limits_{\tau\in U_n(k)} \sgn(\tau)=
\sum\limits_{k=0}^n {n \brack k}_{-1} (-q)^k .
$$
From the case $m=2$ of Corollary~\ref{t.olive} we have 
$$
{n \brack k}_{-1}=
\cases{%
0, &if $k$ and $n-k$ are odd;\cr
\noalign{\vskip4pt}
\displaystyle{\lfloor n/2\rfloor \choose \lfloor k/2\rfloor}, &otherwise.}
$$
Thus, for $n$ even ($n=2m$):
$$
\sum\limits_{\tau\in U_n} \sgn(\tau) q^{\neg(\tau)}=
\sum\limits_{t=0}^m {m\choose t} (-q)^{2t}=(1+q^2)^m,
$$
and for $n$ odd ($n=2m+1$):
$$
\sum\limits_{\tau\in U_n} \sgn(\tau) q^{\neg(\tau)}=
(1-q)\sum\limits_{t=0}^m {m\choose t} (-q)^{2t}=(1-q)(1+q^2)^m .
$$
We conclude that, for $n$ odd ($n=2m+1$):
\begin{eqnarray*}
\sum\limits_{\sigma\in B_n} \sgn(\sigma) q^{\fm(\sigma)}
&=& (1-q)(1+q^2)^m [1]_{q^2}[2]_{-q^2}\cdots [2m+1]_{q^2}\\
&=& (1-q)(1+q^2)^m {\prod\limits_{t=1}^{2m+1}(1-q^{2t})\over
    (1-q^2)^{m+1}(1+q^2)^m}\\
&=& {(1-q)\prod\limits_{t=1}^{2m+1}(1-q^{2t})\over (1-q^2)^{m+1}}
= {\prod\limits_{t=1}^{2m+1}(1-q^{2t})\over (1+q)^{m+1}(1-q)^{m}}\\[3pt]
&=& [2]_{-q} [4]_{q} \cdots [2(2m+1)]_{-q}.
\end{eqnarray*}
The case of $n$ even is similar.

\qed

\section{Other One-Dimensional Characters of $B_n$}\label{s.char}

The group $B_n$ has four one-dimensional characters:
the trivial character; the sign character; $(-1)^{\neg(\sigma)}$;
and the sign of $(|\sigma(1)|,\dots,|\sigma(n)|) \in S_n$,
denoted $\sgn(|\sigma|)$.
We now generalize the results of the previous section to the 
last two one-dimensional characters.

\begin{thm}\label{t.fm2}
$$
\sum\limits_{\sigma\in B_n} (-1)^{\neg(\sigma)} q^{\fm(\sigma)}=
[2]_{-q} [4]_{-q} \cdots [2n]_{- q}.
$$
\end{thm}

\noindent{\bf Proof.} 
Replace $q$ by $-q$ in Theorem \ref{t.ar}, and use the fact that 
the parity of $\fm$ is equal to the parity of $\neg$.

\qed


\begin{thm}\label{t.fm3}
$$
\sum\limits_{\sigma\in B_n} \sgn(|\sigma|) q^{\fm(\sigma)}=
[2]_{q} [4]_{-q} \cdots [2n]_{(-1)^{n-1} q}.
$$
\end{thm}

\noindent{\bf Proof.} Similarly, replace $q$ by $-q$ in Theorem~\ref{t.fm1} and apply the identity 
$\sgn(\sigma)=\sgn(|\sigma|) \cdot (-1)^{\neg(\sigma)}$.

\qed

%

\section{Major Index on Subgroups}\label{s.subgps}

Let $A_n$ be the group of even permutations on $n$ letters.
Then

\begin{cor}\label{t.subgps1}
$$
\sum_{\pi\in A_n}\, q^{\maj(\pi)} =  
{1\over 2}([1]_{q} [2]_{q} \cdots [n]_{q}+
 [1]_{q} [2]_{-q} \cdots [n]_{(-1)^{n-1} q}).
$$
\end{cor}

\noindent{\bf Proof.} 
Clearly,
\begin{eqnarray*}
\sum_{\pi\in A_n} q^{\maj(\pi)}&=&
\sum_{\pi\in S_n}\, {1+\sgn(\pi)\over 2}\ q^{\maj(\pi)}\\
&=&
{1\over 2}\biggl(\sum_{\pi\in S_n}\, q^{\maj(\pi)}
 + \sum_{\pi\in S_n} \sgn(\pi)\, q^{\maj(\pi)}\biggr). 
\end{eqnarray*}
Corollary~\ref{t.maj3} completes the proof.

\qed

\medskip

Let 
$B^+_n$ be the subgroup of even elements in $B_n$,
$D_n$ the subgroup of elements with even $\neg$
(this is a classical Weyl group),
and $C_2\wr A_n$ the subgroup of elements  $\sigma\in B_n$ with even 
$\sgn(|\sigma|)$. Then

\begin{cor}\label{t.subgps2}
$$
\sum_{\sigma\in B^+_n}\, q^{\fm(\sigma)} =  
{1\over 2}([2]_{q} [4]_{q} \cdots [2n]_{q}+
 [2]_{-q} [4]_{q} \cdots [2n]_{(-1)^{n} q}). \leqno(1)
$$
$$
\sum_{\sigma\in D_n}\, q^{\fm(\sigma)} =  
{1\over 2}([2]_{q} [4]_{q} \cdots [2n]_{q}+
 [2]_{-q} [4]_{-q} \cdots [2n]_{-q}). \leqno(2)
$$
$$
\sum_{\sigma\in C_2\wr A_n}\, q^{\fm(\sigma)} =  
{1\over 2}([2]_{q} [4]_{q} \cdots [2n]_{q}+
 [2]_{q} [4]_{-q} \cdots [2n]_{(-1)^{n-1} q}). \leqno(3)
$$
\end{cor}

\noindent{\bf Proof.} Theorem~\ref{t.fm1} implies (1), Theorem~\ref{t.fm2}
implies (2), and Theorem~\ref{t.fm3} implies (3).

\qed

\end{document}